\documentclass[11pt,a4paper]{article}

\usepackage[T1]{fontenc}
\usepackage[utf8]{inputenc}
\usepackage[english]{babel}
\usepackage{lmodern}
\usepackage{amsmath,amssymb,amsfonts,amsthm}
\usepackage{graphicx}
\usepackage{xcolor}
\usepackage{tikz}
\usetikzlibrary{decorations.pathreplacing,automata,calc,positioning,arrows}
\usepackage{enumerate,enumitem,array}
\usepackage[hidelinks]{hyperref}
\usepackage[a4paper,margin=2.7cm]{geometry}

\newtheorem{thm}{Theorem}[section]

\newtheorem{lema}[thm]{Lemma}

\newtheorem{prob}{Problem}

\newtheorem{obs}[thm]{Observation}
\newtheorem{prop}[thm]{Proposition}
\newtheorem{cor}[thm]{Corollary}

\newcommand{\efface}[1]{}

\title{Bootstrap percolation and $P_3$-hull number in direct products of graphs}
\author{Boštjan Bre\v{s}ar\thanks{Faculty of Natural Sciences and Mathematics, University of Maribor, Slovenia; Institute of Mathematics, Physics and Mechanics, Ljubljana, Slovenia. Email: \texttt{bostjan.bresar@um.si}}
\and Jaka Hed\v{z}et\thanks{Institute of Mathematics, Physics and Mechanics, Ljubljana, Slovenia; Faculty of Natural Sciences and Mathematics, University of Maribor, Slovenia. Email: \texttt{jaka.hedzet@imfm.si}}
\and Rebekah Hermann\thanks{Department of Industrial and Systems Engineering, University of Tennessee Knoxville, U.S.A. Email: \texttt{rherrma2@utk.edu}}}
\date{}

\begin{document}
\maketitle

\begin{abstract}
The $r$-neighbor bootstrap percolation is a graph infection process based on the update rule by which a vertex with $r$ infected neighbors becomes infected. We say that an initial set of infected vertices propagates if all vertices of a graph $G$ are eventually infected, and the minimum cardinality of such a set in $G$ is called the $r$-bootstrap percolation number, $m(G,r)$, of $G$. In this paper, we study percolating sets in direct products of graphs. While in general graphs there is no non-trivial upper bound on $m(G\times H,r)$, we prove several upper bounds under the assumption  $\delta(G)\ge r$. 

We also characterize the connected graphs $G$
and $H$ with minimum degree $2$ that satisfy $m(G \times H, 2) = \frac{|V(G \times H)|}{2}$. In addition, we determine the exact values of $m(P_n \times P_m, 2)$, which are $m+n$ if $m$ and $n$ are both even, and $m+n-1$ otherwise.
\end{abstract}



\section{Introduction}

Given a graph $G$ and an integer $r \geq 2$, the {\em $r$-neighbor bootstrap percolation} is a procedure concerned with the states of its vertices which are of two types: {\em infected} or {\em uninfected}.
An initial set of vertices $A_0\ne \emptyset$ is infected, and by the update rule an uninfected vertex becomes infected whenever it has at least $r$ infected neighbors. Once infected, a vertex never changes its state. Given a set $A_0\subseteq V(G)$ of initially infected vertices, if all vertices of $G$ become infected after an update rule is applied sufficiently many times, then $A_0$ is an $r$-{\em percolating set} of $G$, and we also say that $A_0$ {\em propagates}. The $r$-{\em bootstrap percolation number}, $m(G,r)$, of $G$, is the smallest cardinality of an $r$-percolating set of $G$. 

Bootstrap percolation in graphs was considered from different perspectives, which often concerned classes of graphs that have a product-like structure. Early studies of $r$-neighbor bootstrap percolation in graphs were on square grids~\cite{bal-2012,bal-1998} and hypercubes (i.e., the Cartesian products of graphs $K_2$)~\cite{bal-2006}. Recently, those investigations were extended to generalizations of the mentioned structures, notably on the $d$-dimensional square grids~\cite{prz-2020}, and the Cartesian powers of complete graphs~\cite{bid-2021}. One of the basic questions of these investigations is the probability that a set of infected vertices will propagate if each vertex is infected with some probability. Nevertheless, in several of these studies extremal values have also been considered, where the goal is to determine the cardinality of the smallest possible set of infected vertices that propagates. Some other papers only focus on this combinatorial aspect of the smallest set which propagates. For instance, the $r$-neighbor bootstrap percolation number was recently investigated in strong products of graphs~\cite{bh-2023+}.

In the case $r=2$, the $2$-neighbor bootstrap percolation number coincides with the $P_3$-hull number, which arises from the so-called $P_3$-convexity introduced in~\cite{cen-2010}. Efficient algorithms for determining the $P_3$-hull number and related invariants in several classes of graph were found in~\cite{cam-2015}, and in~\cite{coe-2015} inaproximability results were proved for several convexity parameters including the $P_3$-hull number. The $P_3$-hull number was recently determined for caterpillar graphs~\cite{gon-2021}, while exact values or bounds for this number were obtained for Kneser graphs~\cite{gri-2021} and graphs with diameter $2$~\cite{cap-2019}. 
Concerning products of graphs and using the language of bootstrap percolation, the $2$-neighbor bootstrap percolation number was recently determined for Hamming graphs~\cite{bre-2020}, and was also studied in various standard graph products~\cite{coe-2019}.
While the $2$-neighbor bootstrap percolation numbers for strong and lexicographic products of graphs are straightforward to prove,  in general, bounds for the Cartesian product of graphs are more involved; see~\cite{coe-2019}. Moreover, the fourth among the standard graph products (see~\cite{ham-book} for an extensive survey on graph products), the direct product, was not considered in~\cite{coe-2019} at all. 

The definition of the {\em direct product} $G\times H$ of graphs $G$ and $H$ is simple: while $V(G\times H)=V(G)\times V(H)$, vertices $(g,h),(g',h')$ are adjacent in $G\times H$ if $gg'\in E(G)$ and $hh'\in E(H)$. Due to its definition, this graph product is also known as the {\em categorical product} and the {\em tensor product}. When considering graph invariants of graph products, the direct product is often the most challenging one among the four standard graph products.  For instance, the notorious Hedetniemi's conjecture on the chromatic number of the direct product of graphs was posed in 1966, and was intensively studied for half a century (see two surveys~\cite{tar, zhu}), until it was only recently disproved~\cite{shitov}. In this paper, we initiate the study of the $r$-neighbor bootstrap percolation in direct products of graphs.

The paper is organized as follows. In the next section, we provide basic definitions and notation. In Section~\ref{sec:basic}, we show several upper bounds on $m(G\times H,r)$. In particular, we prove that if $G$ is a graph with $\delta(G)\ge r$ and $H$ is a connected graph, then $m(G\times H,r) $ is bounded from above by the order of $G$. Then, in Section~\ref{sec:extremal}, we characterize the connected graphs $G$ with $\delta(G)\ge 2$ for which $m(G\times K_2,2)= |V(G)|$. In Section~\ref{sec:families}, we 
consider the variation of the grid that pertains to direct product, notably graphs  $P_m \times P_n$, for which we obtain closed formulas for the values of their $2$-neighbor bootstrap percolation numbers (and $P_3$-hull numbers). Finally, we conclude with several open problems in Section~\ref{sec:conclusion}. 

\section{Definitions and notation}\label{sec:def}

Let $G$ be a graph and $v\in V(G)$ one of its vertices. By $N_G(v)$, or simply $N(v)$, we denote the {\em neighborhood} of $v$, which is defined as the set of vertices in $V(G)$ that are adjacent to $v$. 
The {\em $r$-neighbor bootstrap percolation} is a color change process that begins with an arbitrarily chosen initial set $A_0 \subseteq V(G)$ of blue vertices, and, for every $t \geq 1$, $A_t = A_{t-1} \cup \lbrace v\in V(G):\, |N(v) \cap A_{t-1}|\geq r \rbrace $ is the set of vertices whose color is blue at time $t$. We say that the set $A_0$ {\em propagates} (or is a {\em percolating set of $G$}) if $\bigcup\limits_{t\ge 0}^{} A_t= V(G)$. 

Given a graph $G$ and $r\geq 2$, let $$m(G,r)=\min \Big\{ |A_0|: \, A_0 \subseteq V(G), \: \bigcup\limits_{t=0}^{\infty}A_t = V(G)\Big\}$$
\noindent be the size of a smallest percolating set in $G$. 
Any percolating set $S$ satisfying $m(G,r) = |S|$ is thus a {\em minimum percolating set}, and $m(G,r)$ is the {\em $r$-percolation number} of $G$. 
We also say that blue vertices are {\em infected}, and so vertices that are not blue are {\em uninfected}, and the process of color change is referred to as the process of infecting vertices in a graph.

It is sometimes convenient to present the percolation procedure in the following alternative way. If $S$ is an $r$-percolating set in $G$, then we may order the vertices of $V(G)\setminus S$ in a sequence $(x_1,\ldots, x_{|V(G)|-|S|})$ such that $x_i$ has at least $r$ neighbors in $S\cup\{x_1,\ldots,x_{i-1}\}$ for all $i\in [|V(G)|-|S|]$.
 
The sequence is not unique, and we may reorder any subsequence of vertices that belong to the same set $A_t$ as defined above, in any order.
 
Let $G$ and $H$ be two graphs, and consider their direct product $G\times H$.
For a fixed $h\in V(H)$, the subset $G^h = \lbrace (g,h): \, g \in V(G) \rbrace$  of $V(G \times H)$ is the {\em $G$-layer} on vertex $h$. Clearly, the subgraph of $G\times H$ induced by the vertices of the $G$-layer on $h$ is the empty graph on $|V(G)|$ vertices. Similarly, for $g \in V(G)$, the {\em $H$-layer} on vertex $g$ is $^g\!H = \lbrace (g,h) : \, h \in V(H) \rbrace$.  It is well known that the direct product $G\times H$ is a connected graph if and only if $G$ and $H$ are both 
connected and at least one of them is non-bipartite~\cite{ham-book}. 

If $G$ is a graph, and $S\subseteq V(G)$, then by $G[S]$ we denote the subgraph of $G$ induced by $S$. For notation that is not explicitly defined here we refer to~\cite{west}. 

%%%%%%%%%%%%%%%%
\section{Upper bounds} 
\label{sec:basic}
%%%%%%%%%%%%%%%

In this section, we present several upper bounds on $m(G\times H,r)$. We start with an observation that such bounds are feasible only if the product of minimum degrees $\delta(G)$ and $\delta(H)$ of $G$ and $H$ is not smaller than $r$. 

For instance, consider the situation when $r=2$ and  $\delta(G) = \delta(H) = 1$.  Note that $\deg_{G\times H}(u,v)=\deg_G(u)\deg_H(v)$, thus $(u,v)$ will have degree $1$ if both $u$ and $v$ have degree $1$ in their graphs. Vertices of degree $1$ are sensitive for the $2$-neighbor bootstrap percolation process, since they cannot become blue by the color change rule, thus they must be colored blue initially (that is, they belong to $A_0$). As an example, let $G$ be the graph $K_{1,n-1}+e$, which is obtained from the star $K_{1,n-1}$ by adding an edge, and let $H\cong K_{1,n-1}$. One can easily see that $(K_{1,n-1}+e)\times K_{1,n-1}$ has $(n-1)(n-3)$ vertices of degree $1$, and with a little more effort that $m\bigl((K_{1,n-1}+e)\times K_{1,n-1},2\bigr)=(n-1)(n-3)+2$, while $|V(G\times H)|=n^2$. Thus, we infer the following observation. 

\begin{obs}
There exists families of graphs $G_n$ and $H_n$, each of order $n$, such that $$\lim_{n\to \infty}{\frac{m(G_n\times H_n,2)}{|V(G_n\times H_n)|}}=1.$$
\end{obs}

\noindent
The above observation readily implies that there exists no constant $c<1$ such that $m(G\times H,2)\le c|V(G\times H)|$ holds for all graphs $G$ and $H$, even if their order is bounded from below.

In the following two upper bounds on $m(G\times H,r)$ we only require that one of the factors (namely, $G$) has minimum degree at least $r$. 

\begin{prop}
\label{prp:basic-upper-bound1}
Let $r\ge 2$ and let $G$ and $H$ be connected graphs. If $\delta(G)\ge r$, then $m(G\times H,r)\le 2m(G,r)$.
\end{prop}
\begin{proof}
Let $h,h'\in V(H)$ be adjacent vertices in $H$, and let $S\subset V(G)$ be a minimum $r$-percolating set in $G$. We claim that $S\times \{h,h'\}$ is an $r$-percolating set in $G\times H$. (Since $|S\times \{h,h'\}|=2m(G,r)$, the claim implies the truth of the proposition.) Let $(x_1,\ldots, x_{|V(G)|-|S|})$  be a sequence of vertices in $V(G)-S$ in an order in which they become blue in an $r$-percolation procedure in $G$ (as noted above, the sequence is not unique, yet we may choose it arbitrarily). Note that in $G\times H$ vertices $(x_i,h)$ and $(x_i,h')$ become blue, since they have at least $r$ blue neighbors in the set $(S\cup\{x_1,\ldots,x_{i-1}\})\times\{h'\}$, respectively $(S\cup\{x_1,\ldots,x_{i-1}\})\times\{h\}$. In this way, all vertices of the layers $G^{h}$ and $G^{h'}$ become blue. Next, since $\delta(G)\ge r$, each vertex of a layer $G^{h''}$, where $h''\in N_H(h)\cup N_H(h')$, has at least $r$  neighbors in $G^{h}\cup G^{h'}$, and so it becomes blue. By induction and connectedness of $H$, we infer that all vertices of $G\times H$ eventually become blue. 
\end{proof}

In the proof of Proposition~\ref{prp:basic-upper-bound1}, the intermediate step of propagation was to have all vertices of a layer colored blue. In the following result, we start with such a set, thus the order of $G$ is an upper bound for $m(G\times H,r)$

\begin{prop}
\label{prp:basic-upper-bound2}
If $G$ is a graph with $\delta(G)\ge r\ge 2$ and $H$ is a connected graph, then $m(G\times H,r)\le |V(G)|$. In particular, any $G$-layer is a percolating set of $G\times H$.
\end{prop}
\begin{proof}
Let $h\in V(H)$ be an arbitrary vertex in $H$, and let $S=G^h$. We claim that $S$ propagates. Indeed, for any $g\in V(G)$ and $h'\in N_H(h)$, vertex $(g,h')$ has at least $r$ neighbors in $G^h$, since $\deg_G(g)\ge  \delta(G)\ge r$. In this way, vertices of $G^{h'}$ become blue at time $t=1$. By using induction and connectedness of $H$, we derive that all vertices of $G\times H$ eventually become blue. Thus, $m(G\times H,r)\le |G^h|=|V(G)|$.
\end{proof}

By exchanging the roles of $G$ and $H$ and using the construction from Proposition~\ref{prp:basic-upper-bound2} for each connected component of $G$ or $H$, we derive the following observation. 

\begin{cor}
\label{cor:upper-bound-r2}
If $G$ and $H$ are graphs with $\min\{\delta(G),\delta(H)\}\ge 2$, then $m(G\times H,2)\le \frac{|V(G\times H)|}{2}.$
\end{cor}

As an example of graphs $G$ and $H$ that attain equality in Corollary~\ref{cor:upper-bound-r2}, we can take $G=C_{2k+1}$ and $H=K_2$. Indeed, $C_{2k+1}\times K_2\cong C_{4k+2}$, and it is easy to see that $m(C_{4k+2},2)=2k+1$.

Note that $\delta(G)\ge r$ and $|V(H)|\ge 3$ imply, by Proposition~\ref{prp:basic-upper-bound2}, that $m(G\times H,r)\le |V(G)|<\frac{|V(G\times H)|}{2}$.
Therefore, in order to characterize the pairs of graphs that attain the bound in Corollary~\ref{cor:upper-bound-r2}, we may restrict to the case when $\delta(G)\ge r$ and $H\cong K_2$. Nevertheless, the solution of this problem leads through the study of graphs that attain the equality in the bound of Proposition~\ref{prp:basic-upper-bound2}, which we consider in the next section.

\section{Graphs with $m(G\times H,2)= |V(G)|$}
\label{sec:extremal}

In this section, we focus on the case $r=2$, and characterize the graphs $G$ and $H$ for which the bound in Proposition~\ref{prp:basic-upper-bound2} is attained.
In order for the equality in the bound of Proposition~\ref{prp:basic-upper-bound2} may have a chance at holding, $G$ must satisfy $\delta(G) = r$. More precisely, the following result holds.
 
\begin{prop}\label{prop:highdegreenbrs} Let $r\ge 2$ and let $G$ and $H$ be connected graphs. If $\delta(G) \geq r$ and there exists a vertex $v\in V(G)$ such that for $r$ of its neighbors $u_1,\ldots u_r$, it holds ${\rm deg}_G(u_i)>r$ for all $i \in [r]$, then $$m(G\times H,r)< |V(G)|.$$
\end{prop}

\begin{proof}
 Let $x$ and $y$ be adjacent vertices of $H$. Color all vertices $(g, x) \in G^x$ blue except for the vertex $(v,x)$, where $v\in V(G)$ has neighbors $u_1,\ldots, u_r$ all of which have degree greater than $r$ in $G$
  (such a vertex $v$ exists by the assumption). This blue set has cardinality $|V(G)|-1$. Now, all vertices $(u_{i}, y)$ are colored blue in the first percolation step since all have at least $r$ distinct blue neighbors of the form $(g,x)$ for some $g\in V(G)$ (since their degree is at least $r+1$). Then $(v,x)$ is colored in the second percolation step since it has $r$ colored neighbors, namely $(u_{1},y),\ldots, (u_{r}, y)$. Therefore all vertices of $G^x$ are blue after the first two steps. By using induction and connectedness of $H$, we derive that all vertices of $G\times H$ become blue. 
 \end{proof}

As a by-product of the above result we have the following observation.

\begin{cor}
    If $G$ is a connected graph with $\delta(G) = d$, then $m(G \times H, r) < |V(G)|$ for all $r < d$. 
\end{cor}

The case for $\Delta(G)=2=\delta(G)$  has already been established (see the observations about cycles). In order to examine other connected graphs with $\Delta(G) \geq 3$, we first need the following auxilliary result, considering connected graphs with $\Delta(G)\geq 3$ of odd order.
 
\begin{lema}\label{lemma:odd number of vertices}
    If $G$ and $H$ are connected graphs with $\delta (G) \ge 2, \Delta (G) \ge 3$ and $|V(G)|=2k+1$, then $m(G\times H,2) \le 2k$. 
\end{lema}

\begin{proof}
   Let $V(H)=\{1,2,\ldots,|V(H)|\}$, where vertices $1$ and $2$ are adjacent. Note that $G$ is not a cycle. 
   
   We claim that $G$ contains a path $P$ such that the induced subgraph $G'=G[V(P)]$ has $\delta(G')\ge 2$ and $G'$ is not an induced cycle. To see that such a path $P$ exists, first note that $G$ has a cycle $C$ as a subgraph. If $C$ is a Hamiltonian cycle, then clearly $C$ is not induced, and for $P$ one can choose a Hamiltonian path. Otherwise, choose an arbitrary $x$ in $C$ that has a neighbor outside $C$ and extend a Hamiltonian path of $C$ with an end-vertex $x$ by traversing along an arbitrary path $P'$ that leaves $C$ from $x$. Since $G$ is finite and $\delta(G)\ge 2$, eventually there will be a vertex from $P'$ that has a neighbor on $P'$ or on $C$, hence the resulting path $P$ has the claimed properties (that is, $G'=G[V(P)]$ has $\delta(G')\ge 2$ and $G'$ is not an induced cycle).  
    
    Now, let $P = v_1v_2\ldots v_p$ be such a path in $G$. Furthermore, let the endvertex $v_1$ be chosen in such a way that either $\deg_{G'}(v_1)>2$ or $\deg_{G'}(v_1)=2$ and $v_1v_p \notin E(G)$. (Indeed, we can do this because if $\deg_{G'}(v_1)=2$ and $v_1v_p \in E(G)$, then $G[V(P)]$ contains a Hamiltonian cycle, which is not induced. Hence, we can choose $v_1$ arbitrarily, as a vertex of degree at least $3$ in $G'$.) 
    By Proposition~\ref{prp:basic-upper-bound2}, it suffices to show that all vertices of a $G$-layer become blue. We consider two cases with respect to the order of $P$. 

    \textbf{Case 1)} $p$ is odd. In this case, color by blue the vertices $$(v_2,1),(v_2,2),(v_4,1),(v_4,2), \ldots ,(v_{p-1},1), (v_{p-1},2).$$ This set enforces that all the vertices in $V(G')\times \{1,2\}$ become blue. Indeed, every vertex $(v_q,i)$ for $3\le q \le p-2$ and $i \in \{1,2\}$ is either initially blue or adjacent to two blue vertices. Note that
    $v_1$ is adjacent to $v_2$, and let $v_1$ be also adjacent to $v_j \neq v_p$. Then vertices $(v_1,1),(v_1,2)$ have two blue neighbors and are colored blue in step $1$. Finally, one can infer that vertices $(v_p,1),(v_p,2)$ also have two blue neighbors and thus become blue.

    Now, let $K=G-G'$. We in addition color by blue all vertices in $V(K) \times \{1\}.$ Once again we get that the whole subgraph induced by $V(G)\times \{1,2\}$ becomes blue, and we infer that the set of vertices initially colored by blue propagates.  Note that the number of these vertices is $2\frac{p-1}{2} + |V(K)| = p-1 + 2k+1 -p = 2k$, which completes the proof of this case. 

    \textbf{Case 2)} $p$ is even. Initially, color by blue all vertices from the set $$V(G') \times \{1\}.$$ Note that all vertices from $V(G')\times \{1,2\}$ become blue due to Proposition~\ref{prp:basic-upper-bound2}.
    Denote $G_0=G'$ and let $G_i$ be obtained from $G_{i-1}$ by adding a maximal path $v_{i_1}v_{i_2}\ldots v_{i_t}$ of vertices in $V(G)\setminus V(G_{i-1})$ such that $v_{i_1}$ is adjacent to a vertex in $G_{i-1}$.     
    (Clearly, unless $G_{i-1}=G$ such a path exists.). Note that $\delta(G_i) \geq 2$ for all $i$. If $t$ is even, then we color blue all vertices $(v_{i_1},1),(v_{i_2},1)\ldots, (v_{i_t},1)$. Once again from Proposition~\ref{prp:basic-upper-bound2} we know that every vertex in $V(G_i) \times \{1,2\}$ becomes blue. Since the order of $G$ is odd this means that for some $i$, there were an odd number of vertices added to $G_{i-1}$ to obtain $G_i$. Without loss of generality let $i=1$, $V(G_0)=\{v_1,v_2,\ldots , v_p\}$, and $V(G_1)=V(G_0) \cup \{v_{p+1},v_{p+2},\ldots ,v_{p+2\ell-1}\}$ for some positive integer $\ell$ and by the construction we may assume $v_{p+j}v_{p+j+1} \in E(G_1)$ for all $1 \leq j \leq 2\ell-2$. Also without loss of generality let $v_pv_{p+1} \in E(G_1)$. Now we color blue the vertices $\{(v_{p+2},1), (v_{p+2},2),$ $(v_{p+4},1), (v_{p+4},2), \ldots ,(v_{p+2\ell-2},1), (v_{p+2\ell-2},2)\}$. Note that vertices $(v_{p+1},1), (v_{p+1},2)$, $(v_{p+3},1), (v_{p+3},2), \ldots , (v_{p+2\ell-3},1), (v_{p+2\ell-3},2)$ have two blue neighbors and are therefore colored blue. Since vertex $v_{p+2\ell-1}$ has degree at least $2$ in $G_1$ this means that vertices $(v_{p+2\ell-1},1), (v_{p+2\ell-1},2)$ are now also adjacent to two blue vertices. This means that we initially colored exactly $p + 2\ell-2 = |V(G_1)|-1$ vertices and all vertices of $V(G_1)\times \{1,2\}$ become blue. 
    Finally, we color blue all vertices in $V(G \setminus G_1) \times \{1\}$. Since the initial set of blue vertices has cardinality $|V(G \setminus G_1)| + |V(G_1)|-1 = |V(G)|-1$ and eventually an entire $G$-layer becomes blue (namely $V(G)\times \{1\}$), this concludes the proof due to Proposition~\ref{prp:basic-upper-bound2}.    
\end{proof}

\begin{lema}\label{lemma:odd cycle}
   If $G$ and $H$ are connected graphs with $\delta (G) = 2, \Delta (G) \ge 3$ such that $G$ contains an odd cycle, then $m(G \times H,2) \leq |V(G)|-1.$ 
\end{lema}

\begin{proof}
    Due to Lemma~\ref{lemma:odd number of vertices}, we only need to consider the case when $G$ has an even order. 
    Let $|V(G)| = n$ and $V(H)=\{1,2,\ldots,|V(H)|\}$, such that $1$ and $2$ are adjacent. In the proof, we will several times use Proposition~\ref{prp:basic-upper-bound2}, without explicitly mentioning it, in the sense that whenever we will notice that vertices of a $G$-layer became blue it yields that the (initially chosen set of blue vertices) propagates.

    Choose a cycle $C$ in $G$, which is isomorphic to $C_p$, where $p$ is odd, and color all vertices of $V(C) \times \{1\}$ blue. Notice that all vertices of $V(C) \times \{2\}$ have two blue neighbors and become blue at step $1$ of the color change process. Let $P^1:v_1v_2\ldots v_t$ be a path in $G$ of maximum length, such that $v_1$ is adjacent to some vertex in $C$. Since the path is maximum, $v_t$ is adjacent to either $v_j$ for some $j\neq t$ or to a vertex in $V(C)$. Now, color by blue all vertices $(v_2,i),(v_4,i), \ldots , (v_t,i)$, if $t$ is even, or $(v_2,i),(v_4,i), \ldots , (v_{t-1},i)$, if $t$ is odd, for $i \in \{1,2\}$. In a similar way as in the proof of Proposition~\ref{lemma:odd number of vertices}, we can show that all vertices of $(V(C) \cup \{v_1,\ldots ,v_t\})\times \{1,2\}$ become blue. 
    If $t$ was even, then we initially colored $p+t$ vertices, and if $t$ was odd, we initially colored $p+t-1$ vertices.     
    Denote by $G_1$ the subgraph of $G$ induced by $V(C)\cup V(P^1)$.
    Notice that if $V(G_1)=V(G)$,  then $t$ must have been odd, so we initially colored $p + 2\cdot \frac{t-1}{2} = p + t-1 = n-1$ vertices and the proof is complete.  Otherwise, we inductively construct a nested sequence of subgraphs $G_i$ of $G$, where $G_1$ has just been defined, and we next define how we obtain $G_{i+1}$  from $G_i$. 

    Thus, let $G_i$ be the subgraph of $G$ such that at most $|V(G_i)|$ vertices of $V(G_i)\times \{1,2\}$ are initially colored blue, and which makes all vertices of $V(G_i)\times \{1,2\}$ be colored blue. Let $P^{i+1}$ be a path of maximum length in $G$ such that a vertex of $P^i$ is adjacent to a vertex in $G_i$. As in the previous paragraph, color  blue the vertices of $V(P^{i+1}) \times \{1,2\}$, starting from the second vertex of $P^{i+1}$ and alternating. Now, the next subgraph in the sequence, $G_{i+1}$, is induced by $V(G_i) \cup V(P^i)$. Note that at most $|V(G_{i+1})|$ vertices have been colored blue initially, and the whole subgraph induced by $V(G_{i+1}) \times \{1,2\}$ became blue after a few steps of applying the color change rule.
    
    Since we started the process with an odd cycle, and $G$ has an even order, 
    it is clear that one of the paths $P^j$ has to have an odd number of vertices. Note that when a path $P^j$ is used in the construction, the graph $G_{j+1}$ defined as induced by $V(G_j) \cup V(P^{j+1})$, and we initially colored at most $|V(G_{j+1})|-1$ vertices. By using the described procedure, one can see that initially at most $|V(G)|-k$ vertices were colored blue, where $k$ is the number of odd paths that appear in the procedure. 
\end{proof}

\begin{figure}[htb]
\begin{center}
\begin{tikzpicture}[scale=1,style=thick,x=1cm,y=1cm]
\def\vr{2.5pt} % \vr = vertex radius;
% define vertices
\path (0,0) coordinate (a);
\path (1,0) coordinate (b);
\path (0.5,0.86) coordinate (c);
\path (2,0) coordinate (d);
\path (3,0) coordinate (e);
\path (2.5,0.86) coordinate (f);

\path (4,0) coordinate (A);
\path (5,0) coordinate (B);
\path (6,0) coordinate (C);
\path (7,0) coordinate (D);
\path (8,0) coordinate (E);
\path (9,0) coordinate (F);
\path (9,1) coordinate (G);
\path (8,1) coordinate (H);
\path (7,1) coordinate (I);
\path (6,1) coordinate (J);
\path (5,1) coordinate (K);
\path (4,1) coordinate (L);

%  edges
\draw (a)--(b)--(c)--(a);
\draw (b)--(d)--(e)--(f)--(d);

\draw (A)--(K)--(C)--(L)--(B)--(J)--(A);
\draw (J)--(D)--(H)--(F)--(I)--(E)--(G)--(D);
\draw (C)--(I);

% vertices
\draw (a) [fill=white] circle (\vr);
\draw (b) [fill=black] circle (\vr);
\draw (c) [fill=black] circle (\vr);
\draw (d) [fill=white] circle (\vr);
\draw (e) [fill=white] circle (\vr);
\draw (f) [fill=black] circle (\vr);

\draw (A) [fill=black] circle (\vr);
\draw (B) [fill=black] circle (\vr);
\draw (C) [fill=black] circle (\vr);
\draw (D) [fill=white] circle (\vr);
\draw (E) [fill=black] circle (\vr);
\draw (F) [fill=white] circle (\vr);
\draw (G) [fill=white] circle (\vr);
\draw (H) [fill=black] circle (\vr);
\draw (I) [fill=white] circle (\vr);
\draw (J) [fill=white] circle (\vr);
\draw (K) [fill=white] circle (\vr);
\draw (L) [fill=white] circle (\vr);

\end{tikzpicture}
\caption{Graphs $G$ and $G\times K_2$ with their percolating sets.}
\label{fig:G and G times K2}
\end{center}
\end{figure}
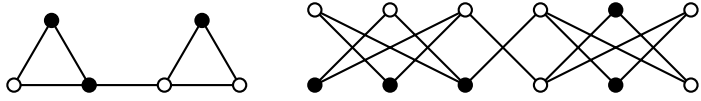

From the previous two lemmas we obtain a result that characterizes the extremal graphs for the bound in Proposition~\ref{prp:basic-upper-bound2} for the case $r=2$. Note that Lemma~\ref{lemma:odd cycle} implies that unless $G$ is an odd cycle, it has to be bipartite. In addition, $G$ has to be a bipartite graph with $m(G,2) \geq \frac{|V(G)|}{2}$. However, since it is clear that in a bipartite graph $G$ with $\delta(G)\ge 2$, we always have $m(G,2)\le \frac{|V(G)|}{2}$, the condition reduces to $m(G,2) =\frac{|V(G)|}{2}$.

\begin{thm}
\label{thm:char}
    If $G$ is a connected graph with $\delta(G) \geq 2$, then $m(G\times K_2,2)=|V(G)|$ if and only if either 
    \begin{itemize}
        \item $G\cong C_{2k+1}$, where $k\ge 2$, or
        \item $G$ is a bipartite graph with $m(G,2) = \frac{|V(G)|}{2}.$
    \end{itemize} 
\end{thm}

As a consequence of the theorem above, we also get a characterization of the graphs $G$ and $H$ with attaining the half of the order bound. 
\begin{cor}
If $G$ and $H$ are connected graphs with $\delta(G) \geq 2$, then \newline $m(G\times H,2)=\frac{|V(G\times H)|}{2}$ if and only if $H\cong K_2$ and either 
    \begin{itemize}
        \item $G\cong C_{2k+1}$, where $k\ge 2$, or
        \item $G$ is a bipartite graph with $m(G,2)=\frac{|V(G)|}{2}.$
    \end{itemize} 
\end{cor}

\section{Direct products of two paths }
\label{sec:families}
%%%%%%%%%%%%%%%%%%%%%%

As mentioned in the introduction and encountered during our studies of upper bounds, 
it is hard to expect any clear general results on the bootstrap percolation numbers of direct product of graphs. Thus, in this section we study a basic instance of direct products, notably, direct products of two paths.

Our aim is to establish the values of $m(P_n\times P_m, 2)$. First, note that whenever $n=2$ (respectively $m=2$) the direct product of two paths is the disjoint union of two paths $P_m$ (respectively $P_n$). Therefore $m(P_2 \times P_m,2)=2 m(P_m,2) = 2\lceil\frac{m}{2}\rceil$. Furthermore, when $n=m=3$, the graph $P_n \times P_m$ is the disjoint union of $C_4$ and $K_{1,4}$, hence $m(P_3 \times P_3,2)=6$. Therefore, in the following results, $n$ and $m$ will always be such that $n\ge m \geq 3$. and not both $m$ and $n$ are equal to $3$. 

In the following lemma, we provide lower bounds for the $2$-bootstrap percolation numbers of direct products of paths. We do this by making use of perimeter arguments that were first used for Cartesian grids by Bollob\'as~\cite{bol-2006}.

\begin{lema}
\label{lem:gridlower}
    Let $n,m \geq 3$ and not both equal to $3$. Then $$m(P_n\times P_m,2)\geq
        \begin{cases}
            n+m; &   \text{both $m$ and $n$ are even   }\\
            n+m-1; &    \text{otherwise}
        \end{cases}$$

\end{lema}

\begin{figure}[htb]
\begin{center}
\begin{tikzpicture}[scale=1,style=thick,x=1cm,y=1cm]
\def\vr{2.5pt} % \vr = vertex radius;
% define vertices
\path (0,4) coordinate (a);
\path (0,1) coordinate (b);
\path (1,0.5) coordinate (c);
\path (1,1.5) coordinate (d);
\path (2,4) coordinate (e);
\path (2,1) coordinate (f);
\path (3,0.5) coordinate (g);
\path (3,1.5) coordinate (h);
\path (4,4) coordinate (i);
\path (4,1) coordinate (j);
\path (5,0.5) coordinate (k);
\path (5,1.5) coordinate (l);
\path (6,4) coordinate (a1);
\path (6,1) coordinate (b1);
\path (7,0.5) coordinate (c1);
\path (7,1.5) coordinate (d1);

\path (0,2) coordinate (a2);
\path (0,3) coordinate (b2);
\path (1,2.5) coordinate (c2);
\path (1,3.5) coordinate (d2);
\path (2,2) coordinate (e2);
\path (2,3) coordinate (f2);
\path (3,2.5) coordinate (g2);
\path (3,3.5) coordinate (h2);
\path (4,2) coordinate (i2);
\path (4,3) coordinate (j2);
\path (5,2.5) coordinate (k2);
\path (5,3.5) coordinate (l2);
\path (6,2) coordinate (a12);
\path (6,3) coordinate (b12);
\path (7,2.5) coordinate (c12);
\path (7,3.5) coordinate (d12);

%  edges
\draw (a)--(d2)--(e)--(h2)--(i)--(l2)--(a1)--(d12);
\draw (b)--(c);
\draw (b)--(d)--(a2);
%\draw (c)--(e);
\draw (c)--(f);
\draw (d)--(f);
%\draw (e)--(g);
\draw (f)--(g);
\draw (f)--(h)--(e2);
%\draw (g)--(i);
\draw (g)--(j);
\draw (h)--(j);
\draw (k)--(b1);
\draw (j)--(k);
\draw (j)--(l)--(b1);
\draw (d1)--(a12)--(l)--(i2);

%\draw (a1)--(c1);
\draw (b1)--(c1);
\draw (b1)--(d1);

\draw (a2)--(c2);
\draw (b2)--(c2);
\draw (b2)--(d2);
\draw (c2)--(e2)--(d);
\draw (c2)--(f2);
\draw (d2)--(f2);
\draw (e2)--(g2);
\draw (f2)--(g2);
\draw (f2)--(h2);
\draw (g2)--(i2)--(h);
\draw (g2)--(j2);
\draw (h2)--(j2);
\draw (i2)--(k2)--(b12);
\draw (j2)--(k2)--(a12);
\draw (j2)--(l2)--(b12);

\draw (a12)--(c12);
\draw (b12)--(c12);
\draw (b12)--(d12);

% vertices
\draw (a) [fill=black] circle (\vr);
\draw (b) [fill=white] circle (\vr);
\draw (c) [fill=black] circle (\vr);
\draw (d) [fill=white] circle (\vr);
\draw (e) [fill=white] circle (\vr);
\draw (f) [fill=white] circle (\vr);
\draw (g) [fill=white] circle (\vr);
\draw (h) [fill=white] circle (\vr);
\draw (i) [fill=white] circle (\vr);
\draw (j) [fill=white] circle (\vr);
\draw (k) [fill=white] circle (\vr);
\draw (l) [fill=white] circle (\vr);
\draw (a1) [fill=white] circle (\vr);
\draw (b1) [fill=white] circle (\vr);
\draw (c1) [fill=black] circle (\vr);
\draw (d1) [fill=white] circle (\vr);

\draw (a2) [fill=black] circle (\vr);
\draw (b2) [fill=white] circle (\vr);
\draw (c2) [fill=white] circle (\vr);
\draw (d2) [fill=white] circle (\vr);
\draw (e2) [fill=white] circle (\vr);
\draw (f2) [fill=black] circle (\vr);
\draw (g2) [fill=white] circle (\vr);
\draw (h2) [fill=white] circle (\vr);
\draw (i2) [fill=white] circle (\vr);
\draw (j2) [fill=black] circle (\vr);
\draw (k2) [fill=white] circle (\vr);
\draw (l2) [fill=white] circle (\vr);
\draw (a12) [fill=black] circle (\vr);
\draw (b12) [fill=white] circle (\vr);
\draw (c12) [fill=white] circle (\vr);
\draw (d12) [fill=black] circle (\vr);

\path (a) node[above] {$(1,1)$};
\path (c1) node[below] {$(8,8)$};
\path (d12) node[above] {$(8,2)$};
\path (c) node[below] {$(2,8)$};

\end{tikzpicture}
\caption{Graph $G_1$ in $P_8\times P_8$ with its percolating set.}
\label{fig:G1 in Pn by Pm}
\end{center}
\end{figure}
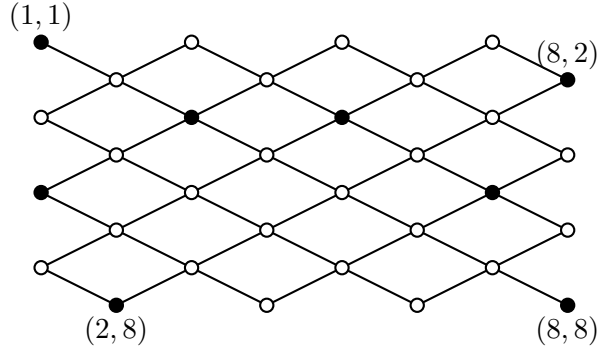

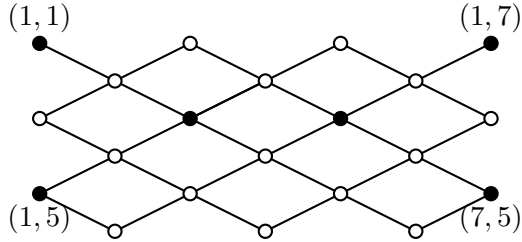
\begin{figure}[htb]
\begin{center}
\begin{tikzpicture}[scale=1,style=thick,x=1cm,y=1cm]
\def\vr{2.5pt} % \vr = vertex radius;
% define vertices
\path (0,1) coordinate (a);
\path (0,0) coordinate (b);
\path (0,-1) coordinate (c);
\path (1,0.5) coordinate (d);
\path (1,-0.5) coordinate (e);
\path (2,1) coordinate (f);
\path (2,0) coordinate (g);
\path (2,-1) coordinate (h);
\path (3,0.5) coordinate (i);
\path (3,-0.5) coordinate (j);
\path (4,1) coordinate (k);
\path (4,0) coordinate (l);
\path (4,-1) coordinate (m);
\path (5,0.5) coordinate (n);
\path (5,-0.5) coordinate (o);
\path (6,1) coordinate (p);
\path (6,0) coordinate (q);
\path (6,-1) coordinate (r);
\path (0,-2) coordinate (s);
\path (1,-1.5) coordinate (t);
\path (2,-2) coordinate (u);
\path (3,-1.5) coordinate (v);
\path (4,-2) coordinate (w);
\path (5,-1.5) coordinate (x);
\path (6,-2) coordinate (y);

%  edges
\draw (a)--(d);
\draw (b)--(d);
\draw (b)--(e);
\draw (c)--(e);
\draw (d)--(f);
\draw (d)--(g);
\draw (e)--(h);
\draw (e)--(i);
\draw (f)--(i);
\draw (g)--(j);
\draw (g)--(k);
\draw (h)--(j);
\draw (i)--(l);
\draw (j)--(l);
\draw (j)--(m);
\draw (k)--(n);
\draw (l)--(n);
\draw (l)--(o);
\draw (m)--(o);
\draw (n)--(p);
\draw (n)--(q);
\draw (o)--(q);
\draw (o)--(r);
\draw (c)--(t);
%\draw (s)--(t);
%\draw (t)--(u);
\draw (t)--(h);
\draw (h)--(v);
%\draw (u)--(v);
%\draw (v)--(w);
\draw (v)--(m);
\draw (m)--(x);
%\draw (w)--(x);
%\draw (x)--(y);
\draw (x)--(r);

% vertices
\draw (a) [fill=black] circle (\vr);
\draw (b) [fill=white] circle (\vr);
\draw (c) [fill=black] circle (\vr);
\draw (d) [fill=white] circle (\vr);
\draw (e) [fill=white] circle (\vr);
\draw (f) [fill=white] circle (\vr);
\draw (g) [fill=black] circle (\vr);
\draw (h) [fill=white] circle (\vr);
\draw (i) [fill=white] circle (\vr);
\draw (j) [fill=white] circle (\vr);
\draw (k) [fill=white] circle (\vr);
\draw (l) [fill=black] circle (\vr);
\draw (m) [fill=white] circle (\vr);
\draw (n) [fill=white] circle (\vr);
\draw (o) [fill=white] circle (\vr);
\draw (p) [fill=black] circle (\vr);
\draw (q) [fill=white] circle (\vr);
\draw (r) [fill=black] circle (\vr);
%\draw (s) [fill=black] circle (\vr);
\draw (t) [fill=white] circle (\vr);
%\draw (u) [fill=white] circle (\vr);
\draw (v) [fill=white] circle (\vr);
%\draw (w) [fill=white] circle (\vr);
\draw (x) [fill=white] circle (\vr);
%\draw (y) [fill=black] circle (\vr);

\path (a) node[above] {$(1,1)$};
\path (p) node[above] {$(1,7)$};
\path (r) node[below] {$(7,5)$};
\path (c) node[below] {$(1,5)$};

\end{tikzpicture}
\caption{Graph $G_2$ in $P_7 \times P_6$ with a percolating set.}
\label{fig:G2 in Pn by Pm}
\end{center}
\end{figure}

\begin{proof}

    Denote $V(P_n)=\{1,2,\ldots,n\}$ and $V(P_m)=\{1,2,\ldots,m\}$. Note that the graph $P_n\times P_m$ is the disjoint union of two graphs. 
    
    When $n$ and $m$ are both even, $P_n\times P_m$ contains two copies of the same graph, which we denote by $G_1$ (see Fig.~\ref{fig:G1 in Pn by Pm} for an example). When $n+m$ is odd, $P_n\times P_m$ contains two copies of the graph, which we denote by $G_2$ (see~Fig.~\ref{fig:G2 in Pn by Pm} for an example). In either case, because of symmetry we will consider the copy of $G_i$ that contains vertex $(1,1)$. (Clearly, $m(P_n\times P_m,2)=2m(G_i,2)$.) Observe that $G_1$ and $G_2$ have a grid-like structure, and we can use the well known perimeter idea of Bollob\'{a}s: we replace the vertices of the graph with unit squares, such that two squares share an edge whenever the corresponding vertices are adjacent. If a vertex is colored blue, then the corresponding square is also blue. Notice that whenever a new square is colored blue by at least two blue neighbors, the total perimeter of the blue domain does not increase, which means that the perimeter of initial blue domain must be at least as big as the total perimeter of the grid. We distinguish three cases.
    
    \textbf{Case 1:} $n,m$ even. Note that only the squares corresponding to vertices $$(1,1),(3,1),\ldots,(n-1,1),(n,2),(n,4),\ldots,(n,m),$$    
$$(n-2,m),(n-4,m), \ldots,(2,m),(1,m-1),(1,m-3),\ldots,(1,3)$$ contribute to the perimeter of the whole grid that arises from $G_1$. Every square contributes $2$ to the total perimeter, except for those corresponding to vertices of degree $1$, which are $(1,1)$ and $(n,m)$, and they contribute $3$. Thus the total perimeter is $3 + 2(n-2) + 2(m-2) + 3 = 2n+2m-2.$ 
    As mentioned earlier, the perimeter of the initial blue domain must also be at least this value. Every square can contribute at most $4$ units to the total perimeter, therefore $m(G,2)=2m(G_1,2) \geq 2 \lceil\frac{2n+2m-2}{4}\rceil = 2 \lceil\frac{n+m-1}{2}\rceil = n+m$.

    \textbf{Case 2:} $n+m$ odd. Without loss of generality, let $n$ be odd. Similarly to Case 1, we count the vertices which contribute to the total perimeter. Those are: 
        $(1,1),(3,1),\ldots,(n,1),(n,3),\ldots,$
        
    \noindent $(n,m-1),(m,n-1),(m,n-3),\ldots,(m,2),(m-1,1),(m-3,1)\ldots,(3,1)$. Each contributes $2$ to the total perimeter except for $(1,1)$ and $(n,1)$, which contribute $3$. The total perimeter is therefore $3+2\lfloor\frac{n-2}{2}\rfloor + 3 + 2(m-2)+2\frac{n-1}{2} = 6+(n-3)+2m-4+(n-1)=2n+2m-2$. We infer that $m(P_n\times P_m,2)=2m(G_1,2) \geq 2 \lceil\frac{2n+2m-2}{4}\rceil = 2 \frac{n+m-1}{2} = n+m-1.$

    \textbf{Case 3:} $n,m$ both odd and $n \geq 5$.
    In this case, the product $P_n \times P_m$ is the disjoint union of two non-isomorphic graphs, $H_1$ and $H_2$. The vertex set of $H_1$ is $\{(2a+1, 2b+1):\,0 \leq a \leq \lfloor \frac{n}{2} \rfloor, 0 \leq b \leq \lfloor \frac{m}{2} \rfloor \} \bigcup \{(2c, 2d):\, 1 \leq c \leq \lfloor \frac{n}{2} \rfloor, 1 \leq d \leq \lfloor \frac{m}{2} \rfloor \}$ 
    See Fig.~\ref{fig:H1 in Pn by Pm} for the example of $H_1$ resulting from $P_7 \times P_5$.
        The vertex set of $H_2$ is $\{(2a+1, 2b ):\,0 \leq a \leq \lfloor \frac{n}{2} \rfloor, 1 \leq b \leq \lfloor \frac{m}{2} \rfloor\} \bigcup \{ (2c, 2d+1):\, 1 \leq c \leq \lfloor \frac{n}{2} \rfloor, 0 \leq d \leq \lfloor \frac{m}{2} \rfloor\}$. Similar to above, vertices $(u,v)$ and $(x,y)$ form an edge if $|u-x| = |v-y| = 1$. See Fig.~\ref{fig:H2 in Pn by Pm}. depicting $H_2$ for $P_7 \times P_5$.

    Again, we shall use the perimeter idea replacing vertices with squares. Let us first count the vertices that contribute to the perimeter of $H_1$. Note there are four squares, corresponding to vertices $(1 ,1 ), (n ,1 ), (1 , m)$ and $(n,m)$, that contribute $3$ to the perimeter of $H_1$, while all other vertices contribute $2$. Counting the number of the edges of the grid of squares that correspond to these vertices, we derive that the perimeter arising from $H_1$ is $12 + 4\lfloor \frac {n-2}{2} \rfloor + 4\lfloor \frac {m-2}{2} \rfloor$.

    Now, in $H_2$, all squares on the perimeter of the corresponding grid contribute exactly $2$ edges to the perimeter. Hence, the total perimeter of $H_2$ is $4\lfloor \frac {n}{2} \rfloor + 4\lfloor \frac {m}{2} \rfloor$. 
  
    Thus, 
    
    \begin{align*}
    m(P_n \times P_m, 2) & \geq \left\lceil \frac{12 + 4\lfloor \frac {n-2}{2} \rfloor + 4\lfloor \frac {m-2}{2} \rfloor}{4} \right\rceil + \left\lceil \frac{4\lfloor \frac {n}{2} \rfloor + 4\lfloor \frac {m}{2} \rfloor}{4} \right\rceil \\
    & \geq \left\lceil \frac{12 + 4\lfloor \frac {n-2}{2} \rfloor + 4\lfloor \frac {m-2}{2} \rfloor}{4} + \frac{4\lfloor \frac {n}{2} \rfloor + 4\lfloor \frac {m}{2} \rfloor}{4}\right\rceil \\
    & = \left\lceil \frac{12 + 2(n-3) + 2(m-3)+ 2(n-1)+ 2(m-1)}{4} \right\rceil \\
    & = \left\lceil \frac{12+2n-6+2m-6+2n-2+2m-2}{4} \right\rceil \\
    &= \left\lceil \frac{4n + 4m -4}{4} \right\rceil \\
    & = n+m-1.
    \end{align*}
\end{proof}

For the upper bounds, we will provide a construction of a 2-percolating set of the appropriate size. 

\begin{lema}
    Let $n,m \geq 3$ and not both equal to $3$. Then $$m(P_n\times P_m,2)\leq
        \begin{cases}
            n+m; &   \text{both $m$ and $n$ are even   }\\
            n+m-1; &    \text{otherwise}
        \end{cases}$$
\end{lema}

\begin{proof}
The proof is by construction of appropriate sets of vertices of the desired cardinality that 2-percolate. We again distinguish three cases: 1) $n$ and $m$ even, 2) $n + m$ odd, and 3) $n,m$ both odd. We use the notation as in Lemma~\ref{lem:gridlower}.

\textbf{Case 1.} $n,m$ even. We will construct a percolating set of $G_1$ and we shall orient
the graph as a grid where $(1, 1)$ is the top left corner, and $(n, m)$ is the bottom right. Consider the following subcases.

\textbf{Subcase 1.1.} $m=4$. Let $A_0=\{(1,1),(n,2),(2,4),(n,4)\} \cup \{(2k+1,3): \, k=1,\ldots,\frac{m-4}{2}\}$. Note that $|A_0|=4+\frac{n-4}{2}=\frac{n+4}{2}$. It can be easily verified that $A_0$ is a percolating set of desired cardinality.   

\textbf{Subcase 1.2.} $m=4k, k>1.$ Let $A_0=\{(1,1),(n,2),(2,m),(n,m)\} \cup \{(2k+1,3) : \, k=1,\ldots,\frac{n-4}{2}\}\cup\{(1,4k+1),(n-1,4k+1) : \, k=1,\ldots,\frac{m-4}{4}\}.$ Then $|A_0|=4+\frac{n-4}{2}+2\frac{m-4}{4}=\frac{n+m}{2}.$

\textbf{Subcase 1.3.} $m=4k+2, k\geq1.$ Let $A_0=\{(1,1),(n,2),(n,m)\} \cup \{(2k+1,3) \, : \, k=1,\ldots,\frac{n-4}{2}\}\cup\{(1,4k+1),(n-1,4k+1) : \, k=1,\ldots,\frac{m-2}{4}\}$. Then $|A_0|=3+\frac{n-4}{2}+2\frac{m-2}{4}=\frac{n+m}{2}.$

In both subcases 1.2 and 1.3, the top four $P_n$-layers, $P_n^1, \ldots, P_n^4$ become infected due to blue vertices in $A_0$ that belong to layers $P_n^1\,\ldots, P_n^5$. In addition to the vertices of the first four $P_n$-layers, also $(1,5),(n-1,5)$ are blue, since they belong to $A_0$. By using these blue vertices, all vertices in $P_n^5$ also become blue. The process continues by gradually infecting vertices of lower levels down to level $8$, so if $m\le 8$, the process is completed. 

Now, for $m> 8$ vertices $(1,9),(n-1,9)$ are also in $A_0$ in both subcases, which ensures that the vertices of layer $P_n^9$ become blue, and the process continues to lower $P_n$-layers. In the bottom two $P_n$-layers we have the following situation. If $m=4k$, then $(2,m)\in A_0$ is needed to resolve vertex $(1,m-1)$, while all other vertices of $P_n^{m-1}$ and $P_n^{m}\setminus\{(n,m),(n-1,m-1),(n-2,m)\}$ become blue due to blue vertices from the level above, while $(n,m)$ is in $A_0$, and $(n-1,m-1)$ and $(n-2,m)$ become blue with a help of vertex $(n,m)$. If $m=4k+2$, the percolation process is straighforward also in the bottom two $P_n$-layers. 

\textbf{Case 2:} $n+m$ odd.

Consider the case when exactly one of $n$ or $m$ is odd, and suppose without loss of generality that $n > m$. As earlier, we shall orient the graph as a grid where $(1,1)$ is the top left corner, $(1,m)$ is the bottom left, $(n,1)$ is the top right, and $(n,m)$ is the bottom right. First note that exactly two of these vertices have degree one and all other vertices have degree 2 or 4. Since the product $P_n \times P_m$ results in two copies of the same graph, $G_2$ (example in Fig.~\ref{fig:G2 in Pn by Pm}), we shall show that $m(G_2, 2) = \frac{n+m-1}{2}$.  %\bb{[[I think in Case 1 we have a different orientation: $(1,1)$ is the bottom left corner, $(1,m)$ top left, $(n,1)$ is bottom right, and $(n,m)$ is the top right. We should unify!]]} 

We distinguish two cases with respect to the parity of $n$.

\textbf{Subcase 2.1: $n$ even.} In this case, the vertices $(1,1)$ and $(1,m)$ have degree one, there are $\lceil \frac{m}{2} \rceil$ vertices in columns that have an odd first coordinate, and $\lfloor \frac{m}{2} \rfloor$ vertices in columns that have an even first coordinate. We shall define a set $B \subset V(G_2)$ and show that it propagates on one copy of $G_2$. Include in $B$ alternating vertices on the line segments $y=x$ and $y = -x + m + 1$ beginning with $x=1$. When $\lceil \frac{m}{2} \rceil$ is odd, end at their intersection $( \lceil \frac{m}{2} \rceil,\lceil \frac{m}{2} \rceil )$.  When $\lceil \frac{m}{2} \rceil$ is even, end at $x = \lfloor \frac{m}{2} \rfloor$, so the two lines do not intersect. Also include in $B$ alternating points from the line segments $y = n+2-x $ and $y = x -n + m -1 $ from $x=n$ until $x = n-\lfloor \frac{m}{2}\rfloor$ when $\lceil \frac{m}{2} \rceil$ is odd (two lines do not intersect at a vertex in $B$) or $x = n - \lfloor \frac{m}{2} \rfloor +1$ when $\lceil \frac{m}{2} \rceil$ is even (two lines intersect at a vertex in $B$). Call these four line segments $\mathcal{L}$. Note the graphs with colored vertices in these two cases ($\lceil \frac{m}{2} \rceil$ odd or even) each have two lines that intersect at a vertex in $B$ and two lines that intersect at a vertex not in $B$. Thus, these two cases are the same up to rotation, so we only consider the first case for the remainder of this case (that is, we assume $\lceil \frac{m}{2} \rceil$ is odd). Add to $B$ all vertices on the line $(j, \lceil \frac{m}{2} \rceil)$ for $ \lceil \frac{m}{2} \rceil + 2 \leq j \leq n-\lceil \frac{m}{2}\rceil$ and $j$ odd, and call this line $\mathcal{L}'$. Thus, the size of $B$ for this subgraph is 
\begin{align*}
    |B| &= \big\lceil \frac{m}{2} \big\rceil + \big\lfloor \frac{m}{2} \big\rfloor + \frac{1}{2}(n - \lceil \frac{ m}{2} \rceil -(\big\lceil \frac{m}{2} \big\rceil + 2) ) + 1\\
    & = m + \frac{1}{2}(n - \frac{m+1}{2}  - \frac{m+1}{2} - 2) + 1 \\
    & = \frac{4m}{4} + \frac{2n - m - 1 - m -1 -4+ 4}{4} \\
    & = \frac{2n - 2m - 2}{4} = \frac{n + m -1}{2}.
 \end{align*}

To see that this set propagates on $G_2$, first note that all vertices in $\mathcal{L}$ that were not initially infected become infected in the first step because they each have two neighbors on $\mathcal{L}$ that are in $B$. 
Let us examine the vertices within the two triangular regions formed by $\mathcal{L}$. In the second percolating step, the vertices $(\lceil \frac{m}{2} \rceil -2, \lceil \frac{m}{2} \rceil)$ and $( n- \lceil \frac{m}{2} \rceil +2,  \lceil \frac{m}{2} \rceil)$ are infected as they each have two neighbors on $\mathcal{L}$. In the $p^{th}$ percolating step (which only occurs if $\frac{m}{2} \geq p$), the vertices $( \lceil \frac{m}{2} \rceil-p-2, y)$ and $( n-\lfloor \frac{m}{2} \rfloor+p+2, y )$ for $ \lceil \frac{m}{2} \rceil - p-2 < y < \lceil \frac{m}{2} \rceil + p+2$ ($y$ odd if $p$ is odd, $y$ even if $p$ is even) are infected. This is because they all have two previously infected neighbors in $( \lfloor \frac{m}{2} \rfloor-p-1, z)$ and $( n-\lfloor \frac{m}{2} \rfloor+p+1,z)$, where $\lceil \frac{m}{2} \rceil - p-1 < z < \lceil \frac{m}{2} \rceil + p+1$, or $\mathcal{L}$, which were previously infected.

To see that vertices outside the triangular regions formed by $\mathcal{L}$ are infected, note the vertices $( j,\lceil \frac{m}{2} \rceil \pm 1)$ for $j \in \{\lceil \frac{m}{2} \rceil + 1, n-\lfloor \frac{m}{2} \rfloor - 1 \}$ are infected in the second step since they have two infected neighbors in $\mathcal{L}' \cup \mathcal{L}$. Once these rows have been infected, it is easy to see that the rows $\lceil \frac{m}{2} \rceil-2$ and $\lceil \frac{m}{2}\rceil +2$ are infected as all vertices in those rows have two neighbors in the rows $\lceil \frac{m}{2} \rceil-1$ and $\lceil \frac{m}{2}\rceil +1$ or one neighbor in those rows and one in $\mathcal{L}$. This process repeats until all rows outside of $\mathcal{L}$ are infected. Since all vertices within the triangular regions formed by $\mathcal{L}$ and all vertices outside of them are eventually infected, $B$ propagates.

Thus, $|A_0| =2|B| = 2\frac{n + m -1}{2} = n+m-1$.

\textbf{Subcase 2.2: $n$ odd.} 
When $n$ is odd, the vertices $(1,1)$ and $(n,1)$ have degree one, and all columns have the same number of vertices. We shall find a set $B \subset V(G_2)$ that propagates on $G_2$. Include in $B$ alternating vertices on the line segments $y=x$ and $y = m-x $ from $x=1$ to $x = \frac{m}{2} -1$ when $\frac{m}{2} $ is even or end at $x = \frac{m}{2}$ when $\frac{m}{2} $ is odd. Also include in $B$ alternating points on the lines $y = n+1-x$ and $y = x -n + m -1 $ from $x=n$ to $x = n-(\frac{m}{2}-2)$ when $\frac{m}{2} $ is even or $x = n-(\frac{m}{2})$ when $\frac{m}{2} $ is odd. When $\frac{m}{2} $ is even, the first pair of lines and second pair of lines do not intersect in $B$, and when $\frac{m}{2} $ is odd, the first pair of lines and second pair of lines intersect at a vertex in $B$. Also include in $B$ all vertices on the line $(j, \frac{m}{2})$ for all feasible $j$, where %\footnote{Before it was written that $j$ is even, but I think this depends on some of the parities of $n$ and $m$. If you do not like "feasible", we could also say: where $j\in V(H_1)$ or something similar.} 
$\frac{m}{2}+2 \leq j \leq n-\frac{m}{2}-1 $. The number of initially infected vertices on this line is $\lfloor \frac{n-m}{2} \rfloor$ since there are $n-m$ columns between the two triangular regions, and only every second $j$ is feasible to appear on this line. In either case, the size of $B$ for this subgraph is 
\begin{align*}
    |B| &= m + \big\lfloor \frac{n-m}{2} \big\rfloor = \frac{n -m -1 + 2m}{2} = \frac{n +m-1}{2}
 \end{align*}

 The percolation process is similar to that in \textbf{Subcase 2.1} with potentially subtle differences near the boundary of $\mathcal{L}$ and $\mathcal{L}'$. Since there are two copies of $G_2$ in $P_n\times P_m$, we get $|A_0| = 2|B| = n+m-1.$

Thus, when $n+m$ is odd, $m(P_n \times P_m, 2) \leq n+m-1$.  

%%%%%%%%%%%%%%%%%%%%%%%%%%%%%%%%%%%%%%%%%%%%%%%%%%%%%%%
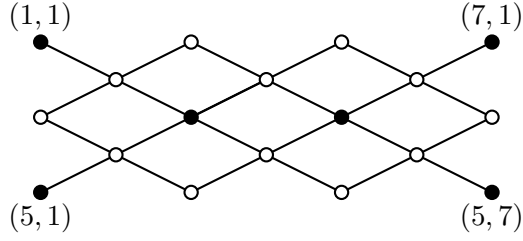
\begin{figure}[htb]
\begin{center}
\begin{tikzpicture}[scale=1,style=thick,x=1cm,y=1cm]
\def\vr{2.5pt} % \vr = vertex radius;
% define vertices
\path (0,1) coordinate (a);
\path (0,0) coordinate (b);
\path (0,-1) coordinate (c);
\path (1,0.5) coordinate (d);
\path (1,-0.5) coordinate (e);
\path (2,1) coordinate (f);
\path (2,0) coordinate (g);
\path (2,-1) coordinate (h);
\path (3,0.5) coordinate (i);
\path (3,-0.5) coordinate (j);
\path (4,1) coordinate (k);
\path (4,0) coordinate (l);
\path (4,-1) coordinate (m);
\path (5,0.5) coordinate (n);
\path (5,-0.5) coordinate (o);
\path (6,1) coordinate (p);
\path (6,0) coordinate (q);
\path (6,-1) coordinate (r);

%  edges
\draw (a)--(d);
\draw (b)--(d);
\draw (b)--(e);
\draw (c)--(e);
\draw (d)--(f);
\draw (d)--(g);
\draw (e)--(h);
\draw (e)--(i);
\draw (f)--(i);
\draw (g)--(j);
\draw (g)--(k);
\draw (h)--(j);
\draw (i)--(l);
\draw (j)--(l);
\draw (j)--(m);
\draw (k)--(n);
\draw (l)--(n);
\draw (l)--(o);
\draw (m)--(o);
\draw (n)--(p);
\draw (n)--(q);
\draw (o)--(q);
\draw (o)--(r);

% vertices
\draw (a) [fill=black] circle (\vr);
\draw (b) [fill=white] circle (\vr);
\draw (c) [fill=black] circle (\vr);
\draw (d) [fill=white] circle (\vr);
\draw (e) [fill=white] circle (\vr);
\draw (f) [fill=white] circle (\vr);
\draw (g) [fill=black] circle (\vr);
\draw (h) [fill=white] circle (\vr);
\draw (i) [fill=white] circle (\vr);
\draw (j) [fill=white] circle (\vr);
\draw (k) [fill=white] circle (\vr);
\draw (l) [fill=black] circle (\vr);
\draw (m) [fill=white] circle (\vr);
\draw (n) [fill=white] circle (\vr);
\draw (o) [fill=white] circle (\vr);
\draw (p) [fill=black] circle (\vr);
\draw (q) [fill=white] circle (\vr);
\draw (r) [fill=black] circle (\vr);

\path (a) node[above] {$(1,1)$};
\path (p) node[above] {$(7,1)$};
\path (c) node[below] {$(5,1)$};
\path (r) node[below] {$(5,7)$};

\end{tikzpicture}
\caption{Graph $H_1$ with a percolating set.}
\label{fig:H1 in Pn by Pm}
\end{center}
\end{figure}

\begin{figure}[htb]
\begin{center}
\begin{tikzpicture}[scale=1,style=thick,x=1cm,y=1cm]
\def\vr{2.5pt} % \vr = vertex radius;
% define vertices
\path (-1,0.5) coordinate (x);
\path (-1,-0.5) coordinate (y);
\path (0,1) coordinate (a);
\path (0,0) coordinate (b);
\path (0,-1) coordinate (c);
\path (1,0.5) coordinate (d);
\path (1,-0.5) coordinate (e);
\path (2,1) coordinate (f);
\path (2,0) coordinate (g);
\path (2,-1) coordinate (h);
\path (3,0.5) coordinate (i);
\path (3,-0.5) coordinate (j);
\path (4,1) coordinate (k);
\path (4,0) coordinate (l);
\path (4,-1) coordinate (m);
\path (5,0.5) coordinate (n);
\path (5,-0.5) coordinate (o);

%  edges
\draw (x)--(a);
\draw (x)--(b);
\draw (y)--(b);
\draw (y)--(c);
\draw (a)--(d);
\draw (b)--(d);
\draw (b)--(e);
\draw (c)--(e);
\draw (d)--(f);
\draw (d)--(g);
\draw (e)--(h);
\draw (e)--(i);
\draw (f)--(i);
\draw (g)--(j);
\draw (g)--(k);
\draw (h)--(j);
\draw (i)--(l);
\draw (j)--(l);
\draw (j)--(m);
\draw (k)--(n);
\draw (l)--(n);
\draw (l)--(o);
\draw (m)--(o);

% vertices
\draw (x) [fill=white] circle (\vr);
\draw (y) [fill=white] circle (\vr);
\draw (a) [fill=black] circle (\vr);
\draw (b) [fill=white] circle (\vr);
\draw (c) [fill=black] circle (\vr);
\draw (d) [fill=white] circle (\vr);
\draw (e) [fill=white] circle (\vr);
\draw (f) [fill=white] circle (\vr);
\draw (g) [fill=black] circle (\vr);
\draw (h) [fill=white] circle (\vr);
\draw (i) [fill=white] circle (\vr);
\draw (j) [fill=white] circle (\vr);
\draw (k) [fill=black] circle (\vr);
\draw (l) [fill=white] circle (\vr);
\draw (m) [fill=black] circle (\vr);
\draw (n) [fill=white] circle (\vr);
\draw (o) [fill=white] circle (\vr);

\path (x) node[above] {$(1,2)$};
\path (y) node[below] {$(1,4)$};
\path (n) node[above] {$(7,2)$};
\path (o) node[below] {$(7,4)$};

\end{tikzpicture}
\caption{Graph $H_2$ with a percolating set.}
\label{fig:H2 in Pn by Pm}
\end{center}
\end{figure}
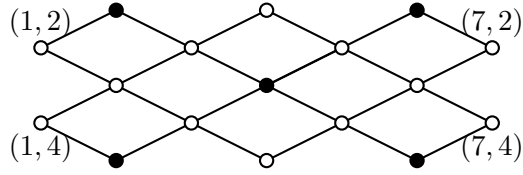

\textbf{Case 3:} $n,m$ both odd, $n \geq 5$. $P_n \times P_m$ results in two disjoint graphs that are not isomorphic. Call them $H_1$ and $H_2$. An example of each subgraph is shown in Figs.~\ref{fig:H1 in Pn by Pm}-\ref{fig:H2 in Pn by Pm}. There are four distinct combinations of $\lfloor \frac{n}{2} \rfloor$ and $\lfloor \frac{m}{2} \rfloor$ even and odd. Some of these combinations have identical 2-percolating sets. We will consider the following three cases:
\begin{enumerate}%[label=3.]
    \item $\lfloor \frac{n}{2} \rfloor$ and $\lfloor \frac{m}{2} \rfloor$ both odd and $n \neq m$ and the equivalent case $\lfloor \frac{n}{2} \rfloor$ even and $\lfloor \frac{m}{2} \rfloor$ odd;
    \item $\lfloor \frac{n}{2} \rfloor$ and $\lfloor \frac{m}{2} \rfloor$ both odd and $n = m$; 
    \item $\lfloor \frac{n}{2} \rfloor$ and $\lfloor \frac{m}{2} \rfloor$ both even and the equivalent case $\lfloor \frac{n}{2} \rfloor$ odd and $\lfloor \frac{m}{2} \rfloor$ even.
\end{enumerate}

\textbf{Subcase 3.1: $\lfloor \frac{n}{2} \rfloor$ and $\lfloor \frac{m}{2} \rfloor$ both odd and $n \neq m$ and the equivalent case $\lfloor \frac{n}{2} \rfloor$ even and $\lfloor \frac{m}{2} \rfloor$ odd.}

We will construct a set $B_0$ in $H_1$ that propagates. Add to $B_0$ alternating vertices on the lines $y = x$ and $y = -x + m+1$ from $x=1$ to $x = \lfloor \frac{m}{2} \rfloor$ and $y=-x+n+1$ and $y=x-n+m$ from $x=n$ until $x = n - \lfloor \frac{m}{2} \rfloor$. Call this set of lines $\mathcal{L}$. In this case, none of the lines intersect within the domain. Then add to $B_0$ vertices $( j,\lceil \frac{m}{2} \rceil)$ for $\lceil \frac{m}{2} \rceil + 2 \leq j \leq n-\lceil \frac{m}{2} \rceil - 2$, and call this set $\mathcal{L}'$. There are $n - 2 \lceil \frac{m}{2} \rceil = n-m-1$ columns that are not in the domain of the four line segments described above. Half of these are contained in $B_0$, rounded down. Thus, the size of $B_0$ is
\begin{align*}
    |B_0| &= 4(\frac{\lceil \frac{m}{2} \rceil}{2}) + \lceil \frac{n-m-1}{2} \rceil \\
    & = 2\frac{m+1}{2}  + \frac{n -m-2}{2}\\
    & = \frac{n+m}{2} = \lceil \frac{n+m-1}{2} \rceil.
\end{align*}

%\rh{finish proof that this 2-propagates $H_1$ as other cases above.}
To see that $B_0$ propagates on $H_1$, first note that all vertices in $\mathcal{L}$ that were not infected initially become infected in the first step because they each have two neighbors on $\mathcal{L}$ that are in $B_0$. The vertices $(\lfloor \frac{m}{2} \rfloor-1 , \lceil \frac{m}{2}\rceil )$ and $(n-\lfloor \frac{m}{2} \rfloor-1 , \lceil \frac{m}{2}\rceil )$ are infected in the first step since the former has neighbors $(\lfloor \frac{m}{2} \rfloor,\lfloor \frac{m}{2} \rfloor), (\lfloor \frac{m}{2} \rfloor , \lceil \frac{m}{2} \rceil + 1) $ and the latter has neighbors $(n-\lfloor \frac{m}{2} \rfloor,\lfloor \frac{m}{2} \rfloor), (n-\lfloor \frac{m}{2} \rfloor , \lceil \frac{m}{2} \rceil + 1) $ that are the last points of the line segments above in $B_0$. Also, all vertices of the form $( k,\lceil \frac{m}{2} \rceil \pm 1)$ for odd $k$ in $\lfloor \frac{m}{2} \rfloor + 3 \leq k \leq n-\lfloor \frac{m}{2} \rfloor - 3$ are infected in the first step since they have two neighbors in $( j,\lceil \frac{m}{2} \rceil) \subset B_0$ for even $j$ in $\lceil \frac{m}{2} \rceil + 3 \leq j \leq n-\lceil \frac{m}{2} \rceil - 3$. 

Let us examine the vertices within each triangular region formed by $\mathcal{L}$. In the second percolating step, the vertices $( \lfloor \frac{m}{2} \rfloor-2, \lceil \frac{m}{2} \rceil \pm 1)$ and $( n-\lfloor \frac{m}{2} \rfloor+2, \lceil \frac{m}{2} \rceil \pm 1 )$ for are infected as they each have one neighbor on $\mathcal{L}$ and one neighbor in $\{(\lfloor \frac{m}{2} \rfloor-1 , \lceil \frac{m}{2}\rceil ), (n-\lfloor \frac{m}{2} \rfloor-1 , \lceil \frac{m}{2}\rceil )\}$.  In the $p^{th}$ percolating step (which only occurs if $\frac{m}{2} \geq p$), the vertices $( \lfloor \frac{m}{2} \rfloor-p, y)$ and $( n-\lfloor \frac{m}{2} \rfloor+p, y )$ for $ \lceil \frac{m}{2} \rceil - p < y < \lceil \frac{m}{2} \rceil + p$ ($y$ odd if $p$ is odd, $y$ even if $p$ is even) are infected. This is because they all have two previously infected neighbors in $(\lfloor \frac{m}{2} \rfloor-p+1, z)$ and $( n-\lfloor \frac{m}{2} \rfloor+p-1,z)$, where $\lceil \frac{m}{2} \rceil - p+1 < z < \lceil \frac{m}{2} \rceil + p-1$, or $\mathcal{L}$, which were previously infected.

To see that vertices outside $\mathcal{L}$ are infected, note the vertices $( j,\lceil \frac{m}{2} \rceil \pm 1)$ for $j \in \{\lfloor \frac{m}{2} \rfloor + 2, n-\lfloor \frac{m}{2} \rfloor - 2 \}$ are infected in the second step since they have two infected neighbors in $\mathcal{L}' \cup \{(\lfloor \frac{m}{2} \rfloor,\lfloor \frac{m}{2} \rfloor), (\lfloor \frac{m}{2} \rfloor , \lceil \frac{m}{2} \rceil + 1),(n-\lfloor \frac{m}{2} \rfloor,\lfloor \frac{m}{2} \rfloor) , (n-\lfloor \frac{m}{2} \rfloor , \lceil \frac{m}{2} \rceil + 1)\}$. In the third step, the vertices $(\lceil \frac{m}{2}\rceil + 1,\lceil \frac{m}{2} \rceil \pm 1 )$ and $(n - \lceil \frac{m}{2}\rceil - 1 ,\lceil \frac{m}{2} \rceil \pm 1)$ are infected since they have neighbors that are in $\mathcal{L}' \cup \{(\lceil \frac{m}{2}\rceil , \lceil \frac{m}{2}\rceil),(n-\lceil \frac{m}{2}\rceil, \lceil \frac{m}{2}\rceil)\}$. Note this means that the rows $\lceil \frac{m}{2} \rceil-1$ and $\lceil \frac{m}{2}\rceil +1$ above $\mathcal{L}$ are infected. Since these rows are infected, the rows $\lceil \frac{m}{2} \rceil-2$ and $\lceil \frac{m}{2}\rceil +2$ are infected as all vertices in those rows have two neighbors in the rows $\lceil \frac{m}{2} \rceil-1$ and $\lceil \frac{m}{2}\rceil +1$ or one neighbor in those rows and one in $\mathcal{L}$. This process repeats until all rows outside of $\mathcal{L}$ are infected. Since all vertices within the triangular regions formed by $\mathcal{L}$ and all vertices outside of them are eventually infected, $B_0$ propagates.

Now we will construct a set $B_1$ that propagates on $H_2$. Add to $B_1$ alternating vertices on the lines $y = x+1 $ and $y = m-x$ from $x=1$ to $x = \lfloor \frac{m}{2} \rfloor$. These two lines intersect at $x = \lfloor \frac{m}{2} \rfloor$. Also add to $B_1$ alternating vertices on the lines $y = -x + n + 2$ and $ y = x -n+m-1 $ from $x = n$ to $x = n-\lfloor \frac{m}{2} \rfloor$. These lines intersect at $x = n-\lfloor \frac{m}{2} \rfloor$ in this case. Furthermore, add to $B_1$ vertices on the line $(j, \lceil \frac{m}{2} \rceil)$ from $ \lfloor \frac{m}{2} \rfloor + 2 \leq j \leq n - \lfloor \frac{m}{2} \rfloor -1 $. There are $n - 2 \lfloor \frac{m}{2} \rfloor = n-m+1$ columns that are not in the domain of the four line segments described above. Half of these are contained in $B_0$, rounded down. The size of $B_1$ is

\begin{align*}
    |B_1| &= 4(\frac{\lfloor \frac{m}{2} \rfloor}{2}) + \lfloor \frac{n-m+1}{2} \rfloor \\
    & = 2\frac{m-1}{2}  + \frac{n -m}{2}\\
    & = \frac{n+m-2}{2} = \lfloor \frac{n+m-1}{2} \rfloor.
\end{align*}

The proof that $B_1$ propagates in $H_2$ can be done in similar manner as the one for $B_0$ on $G_1$, by taking into account that the values of $j$ for the vertices in the second percolating step are $j \in \{ \lfloor \frac{m}{2} \rfloor + 1, n - \lfloor \frac{m}{2} \rfloor -1\}.$

Thus, $B_0 \cup B_1$ is a 2-percolating set in $P_n \times P_m$ and $|B_0| + |B_1| = n+m-1$.

\textbf{Subcase 3.2: $\lfloor \frac{n}{2} \rfloor$ and $\lfloor \frac{m}{2} \rfloor$ both odd and $n = m$.}

The case $n=m$ requires a slight variation in either $B_0$ or $B_1$. This variation is required because when $n=m$, the set $B_0$ as defined in~\textbf{Subcase 3.1} has $m+1 = n+1$ vertices, while $B_1$ still has $m-1$ vertices. This gives a total of $n+m$ vertices. We shall provide an alternative construction for $B_0$ here that has $n$ vertices and 2-propagates. 

Add to $B_0$ alternating vertices on the lines $y = x$ and $y = -x + m+1$ from $x=1$ to $x = \lfloor \frac{m}{2} \rfloor - 2$, and $y=-x+n+1$ and $y=x-n+m$ from $x=n$ until $x = n - \lfloor \frac{m}{2} \rfloor + 2$. In this case, none of the lines intersect within the domain. Then add to $B_0$ the set of vertices in $A = \{( \lceil \frac{m}{2} \rceil, j):\, \lfloor \frac{m}{2} \rfloor-2 \leq j \leq \lfloor \frac{m}{2} \rfloor +2 \}$  that are in $V(H_1)$. Note $|V(H_1) \cap A| = 3$. The size of $B_0$ is

\begin{align*}
    |B_0| & = 4(\frac{\big\lceil \frac{m}{2} \big\rceil - 2}{2} ) + 3 \\
    & = 2(\frac{m+1}{2}-2) + 3 \\
    & = m = n.
\end{align*}

This set propagates similarly to \textbf{Subcase 3.1}.
Thus, $A_0 = B_0 \cup B_1$ propagates in $P_n\times P_m$, and $|A_0| = |B_0| + |B_1| = n + m-1$

\textbf{Subcase 3.3: $\lfloor \frac{n}{2} \rfloor$ and $\lfloor \frac{m}{2} \rfloor$ both even and the equivalent case $\lfloor \frac{n}{2} \rfloor$ odd and $\lfloor \frac{m}{2} \rfloor$ even}. %\rh{this subcase is finished except for the 2-percolation proofs}

%\rh{odd $\lfloor \frac{n}{2} \rfloor$ odd and  $\lfloor \frac{m}{2} \rfloor$ even have a clear pattern and $m=5$ upholds it, so we can remove that case and put into this one. }

Let us construct the set of vertices $B_0$ that propagates in $H_1$. Add to $B_0$ alternating vertices on the lines $y = x$ and $y = -x + m+1$ from $x=1$ to $x = \lceil \frac{m}{2}\rceil$, and $y=-x+n+1$ and $y=x-n+m$ from $x=n$ until $x = n - \lfloor \frac{m}{2}\rfloor$. The former two line segments contain a total of $\lfloor\frac{m}{2} \rfloor + 1 =\lceil \frac{m}{2}\rceil$ vertices of $B_0$, and the latter two line segments contain a total of $\lfloor\frac{m}{2} \rfloor + 1 = \lceil \frac{m}{2}\rceil$ vertices of $B_0$. Thus, the number of vertices in these four line segments and $B_0$ is $2\lceil \frac{m}{2}\rceil$. Then add to $B_0$ vertices $( j,\lceil \frac{m}{2} \rceil)$ for $\lceil \frac{m}{2}\rceil+2 \leq j \leq n-\lfloor \frac{m}{2}\rfloor- 2$. There are $n - 2 \lceil \frac{m}{2} \rceil = n-m-1$ columns that are not in the domain of the four line segments described above. Half of these are contained in $B_0$, rounded down. 
The size of $B_0$ is thus
\begin{align*}
    |B_0| & = 2\lceil \frac{m}{2} \rceil + \lfloor \frac{n-m-1}{2} \rfloor \\
    & = 2\lceil \frac{m}{2} \rceil + \frac{n-m-2}{2} \\
    &= \frac{2(m+1)+ n - m -2}{2} \\
    &= \frac{n+m}{2} = \big\lceil\frac{n+m-1}{2} \big\rceil
\end{align*}

Now let us define $B_1$ and show that it propagates in $H_2$. Add to $B_1$ alternating vertices on the lines $y = x+1 $ and $y = m-x$ from $x=1$ to $x = \lfloor \frac{m}{2} \rfloor$. Also add to $B_1$ alternating vertices on the lines $y = -x + n + 2$ and $ y = x -n+m-1 $ from $x = n$ to $x = n-\lfloor \frac{m}{2} \rfloor$. Finally, add vertices on the line $( j, \lceil \frac{m}{2} \rceil)$ from $ \lfloor \frac{m}{2} \rfloor +2 \leq j \leq n - \lfloor \frac{m}{2} \rfloor -2 $. There are $n - 2 \lfloor \frac{m}{2} \rfloor = n-m+1$ columns that are not in the domain of the four line segments described above that contain vertices of $B_1$. Half of these are contained in $B_0$, rounded down.  The size of $B_1$ is then
\begin{align*}
    |B_1| & = 4 \frac{\lfloor \frac{m}{2}\rfloor}{2} + \lfloor \frac{n-m+1}{2}\rfloor \\
    & = 2\lfloor \frac{m}{2} \rfloor + \frac{n-m}{2} \\
    &= \frac{2(m-1)+ n - m}{2} \\
    &= \frac{n+m-2}{2} = \big\lfloor\frac{n+m-1}{2} \big\rfloor
\end{align*}

The sets $B_0$ and $B_1$ propagate in a similar was as in \textbf{Subcase 3.1}. %\rh{finish proof that this 2-propagates $H_2$ as other cases above.}
Thus, $B_0 \cup B_1$ is a percolating set of $P_n \times P_m$ and $|B_0| + |B_1| = n+m-1$.

    \qed
\end{proof}

From the above lemmas and the fact that $m(P_3\times P_3,2)=6$, we derive the following result. 

\begin{thm}
\label{thm:grid}
    Let $n,m \geq 3$ and not both equal to $3$. Then $$m(P_n\times P_m,2) =
        \begin{cases}
            n+m; &   \text{both $m$ and $n$ are even   }\\
            n+m-1; &    \text{otherwise}
        \end{cases}$$

\end{thm}

\section{Concluding remarks}\label{sec:conclusion}

\begin{figure}[htb]
\begin{center}
\begin{tikzpicture}[scale=0.9,style=thick,x=1cm,y=1cm]
\def\vr{2.5pt} % \vr = vertex radius;
% define vertices
\path (0,0) coordinate (a);
\path (2,1.5) coordinate (b);
\path (2,0) coordinate (c);
\path (2,-1.5) coordinate (d);
\path (4,1.5) coordinate (e);
\path (4,0) coordinate (f);
\path (4,-1.5) coordinate (h);
\path (6,0) coordinate (i);

\path (8,0) coordinate (A);
\path (10,1.8) coordinate (B);
\path (10,0.8) coordinate (C);
\path (10,-0.8) coordinate (D);
\path (10,-1.8) coordinate (E);
\path (12,1.8) coordinate (F);
\path (12,0.8) coordinate (G);
\path (12,-0.8) coordinate (I);
\path (12,-1.8) coordinate (J);
\path (14,0) coordinate (K);

%  edges
\draw (i)--(e)--(b)--(a)--(d)--(h)--(i);
\draw (a)--(c)--(f)--(i);

\draw (A)--(B)--(F)--(K)--(G)--(C)--(A);
 \draw (A)--(D)--(I)--(K)--(J)--(E)--(A);

% vertices
\draw (a) [fill=white] circle (\vr);
\draw (b) [fill=black] circle (\vr);
\draw (c) [fill=black] circle (\vr);
\draw (d) [fill=black] circle (\vr);
\draw (e) [fill=white] circle (\vr);
\draw (f) [fill=white] circle (\vr);
\draw (h) [fill=white] circle (\vr);
\draw (i) [fill=black] circle (\vr);

\draw (A) [fill=white] circle (\vr);
\draw (B) [fill=black] circle (\vr);
\draw (C) [fill=black] circle (\vr);
\draw (D) [fill=white] circle (\vr);
\draw (E) [fill=white] circle (\vr);
\draw (F) [fill=white] circle (\vr);
\draw (G) [fill=white] circle (\vr);
\draw (I) [fill=black] circle (\vr);
\draw (J) [fill=black] circle (\vr);
\draw (K) [fill=white] circle (\vr);

     \end{tikzpicture}
\caption{Bipartite graphs $B_2$ and $B_3$, where $m(B_2,2)=|V(B_2)|/2$ and $m(B_3,2)<|V(B_3)|/2$.}
\label{fig:bipartite}
\end{center}
\end{figure}
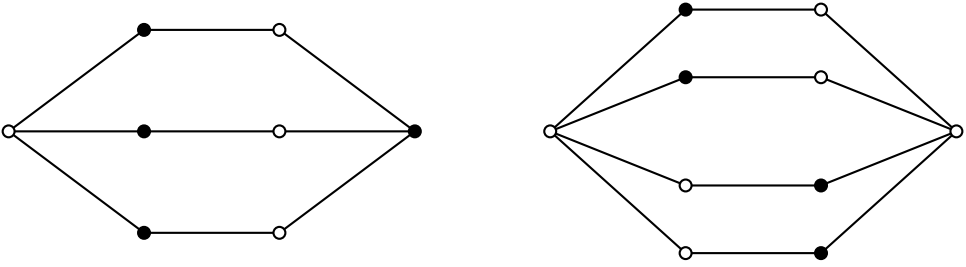

Several open problems naturally arise from the study in this paper. In Theorem~\ref{thm:char} we characterized the graphs $G$ and $H$ for which 
$m(G\times K_2,2)=|V(G)|$. In the characterization appears the class of 
bipartite graphs, which satisfy $m(G,2)=\frac{|V(G)|}{2}$, whose structure has yet to be discovered.
In Fig.~\ref{fig:bipartite} one can find two examples of bipartite graphs $B_2$ and $B_3$, which belong to an infinite family starting with $B_1\cong C_6$. $B_m$ is defined to be $C_6$ plus $m-1$ disjoint paths of length three between exactly one pair of antipodal vertices. Among all graphs from the family, we have $m(B_n,2)=|V(B_n)|/2$ only for $n\in \{1,2\}$. This indicates that the characterization of bipartite graphs that appear in Theorem~\ref{thm:char} might be a challenging problem, which we next propose.

\begin{prob}
Characterize bipartite graphs $G$ with $\delta(G)\ge 2$ satisfying $m(G,2)=\frac{|V(G)|}{2}$.
\end{prob}

It would also be interesting to extend the consideration of graphs with $m(G\times H,r)=|V(G)|$ to $r>2$. 
 A simple family of graphs that satisfy this equality with an arbitrarily large $r$ is that of uniform complete bipartite graphs, $K_{r,r}$, multiplied with $K_2$. Note that $m(K_{r,r},r) = r$ and that $K_{r,r} \times K_2$ consists of two disjoint copies of $K_{r,r}$. Thus, $m(K_{r,r} \times K_2, r) = 2m(K_{r,r},r) = 2r = |V(K_{r,r})|$. 

Another interesting problem is to determine which graphs achieve equality in Proposition~\ref{prp:basic-upper-bound1}.

\begin{prob}
\label{prob:mGxH=2m} 
For which connected graphs $G$ and $H$, where $\delta(G)\ge r\ge 2$,  we have $$m(G\times H,r)=2m(G,r)?$$
\end{prob}

Consider the family of hypercubes $Q_n$. Note that $m(Q_n, n) = 2^{n-1}$ for all $n \geq 2$, as $m(Q_n, r) \geq 2^{r-1}$ (\cite[Theorem 1.3]{morrison}), and equality can be easily found by coloring exactly half of the vertices of $Q_n$ that form a maximum independent set. Since $Q_n \times K_2$ is two disjoint copies of $Q_n$, clearly $m (Q_n \times K_2) = 2\cdot2^{n-1}$, hence hypercubes satisfy the equality in Problem~\ref{prob:mGxH=2m}.

It is quite common to look for bounds on graph invariants of graph products that are expressed in terms of the same invariants of their factors. In this line of study, we call a graph invariant $\sigma$ {\em submultiplicative} (resp. {\em supermultiplicative}) if for every graphs $G$ and $H$ we have $\sigma(G\times H)\le \sigma(G)\sigma(H)$ (resp. $\sigma(G\times H)\ge \sigma(G)\sigma(H)$). We suspect that for graphs $G$ and $H$ (perhaps under some additional boundary restrictions) the $2$-percolation number is submultiplicative, which we pose as the following problem.

\begin{prob}
\label{prob:productbound}
Let $G,H$ be connected graphs, one of which is not bipartite. Is it true that $$m(G\times H, r) \le m(G,r) m(H,r)?$$
\end{prob}

From the example of $H=P_3$ and $G$ is the cycle $C_3$ with a pendant leaf we can see, that the set $S_G \times S_H$, where $S_G, S_H$ are the percolating sets of $G, H$ respectively, is in general not a percolating set of $G \times H$. However $m(G \times H,2)=m(G,2)m(H,2)$ still holds in this example. 

We also ask the following question, which intuitively should hold true, but we have not found a formal proof.

\begin{prob}
Is it true that $\min \{m(G, r), m(H,r)\} \leq m(G \times H, r)$?
\end{prob}

In Section~\ref{sec:families} we considered the $2$-bootstrap percolation numbers of products of paths. The following problem is the next step in this line of investigation.  
\begin{prob}\label{prob:cyclepath}
Determine $m(C_n \times P_m,2)$.  
\end{prob}
From Proposition~\ref{prp:basic-upper-bound2} we see that $m(C_n \times P_m,2)\leq n$, and equality appears to hold for small examples. In fact, when $m=2$, equality holds since when $n$ is odd, $C_n \times P_2 = C_{2n}$ and when $n$ is even, $C_n \times P_2$ is two disjoint copies of $C_n$. Thus, the problem need only be considered for $n \geq 3, m \geq 3$. Establishing a general lower bound for $ m(G \times H, r)$ could potentially solve this problem. 

We were also considering direct products with complete graphs, where our preliminary result states that $m(G\times K_n,r)=r$ whenever $G$ is a connected graph of order at least $3$ and $n \geq 2r$. In particular, this implies that $m(G\times K_n,2)=2$ for any connected graph $G$ and $n\ge 4$. In this direction, we are left to determine the $2$-bootstrap percolation numbers of graphs $G\times K_3$.

A majority of results in this paper are concerned with the $r$-neighbor bootstrap percolation number where $r=2$, which coincides with the $P_3$-hull number. In particular, Theorem~\ref{thm:char} characterizes connected graphs $G$ with $\delta(G) \ge 2$ for which the $P_3$-hull number equals the order of $G$, while Theorem~\ref{thm:grid} provides the values of the $P_3$-hull numbers in direct products of paths. 

\section*{Acknowledgements}
B.B. and J.H. were supported by the Slovenian Research and Innovation Agency (ARIS) under the grants P1-0297, N1-0285, and J1-4008.


\begin{thebibliography}{99}

\bibitem{bal-2006} J.~Balogh, B.~Bollob\'{a}s. Bootstrap percolation on the hypercube. \emph{Probab.\ Theory Related Fields} \textbf{134} (2006), 624--648. 
\href{https://doi.org/10.1007/s00440-005-0451-6}{doi:10.1007/s00440-005-0451-6}

\bibitem{bal-2012} J.~Balogh, B.~Bollob\'{a}s, H.~Duminil-Copin, R.~Morris. The sharp threshold for bootstrap percolation in all dimensions. \emph{Trans.\ Amer.\ Math.\ Soc.} \textbf{364} (2012), 2667--2701.\\ \href{https://doi.org/10.48550/arXiv.1010.3326}{doi:10.48550/arXiv.1010.3326}

\bibitem{bal-1998} J.~Balogh, G.~Pete. Random disease on the square grid. \emph{Random Structures Algorithms} \textbf{13} (1998), 409--422.

\bibitem{bid-2021} M.~Bidgoli, A.~Mohammadian, B.~Tayfeh-Rezaie. Percolating sets in bootstrap percolation on the Hamming graphs and triangular graphs. \emph{European J.\ Combin.} \textbf{92} (2021), \#103256.
\href{https://doi.org/10.1016/j.ejc.2020.103256}{doi:10.1016/j.ejc.2020.103256}



\bibitem{bol-2006} B.~Bollob\'{a}s. \emph{The Art of Mathematics: Coffee Time in Memphis}. Cambridge University Press, 2006.

\bibitem{bh-2023+} B.~Bre\v sar, J.~Hed\v{z}et. Bootstrap percolation in strong products of graphs. \emph{Electron.\ J.\ Combin.} \textbf{31} (2024), \#P4.35.
\href{https://doi.org/10.37236/11826}{doi:10.37236/11826}


\bibitem{bre-2020} B.~Bre\v sar, M.~Valencia-Pabon. On the $P_3$-hull number of Hamming graphs. \emph{Discrete Appl.\ Math.} \textbf{282} (2020), 48--52.
\href{https://doi.org/10.1016/j.dam.2019.11.011}{doi:10.1016/j.dam.2019.11.011}

\bibitem{cam-2015} V.~Campos, R.M.~Sampaio, A.~Silva, J.L.~Szwarcfiter. Graphs with few $P_4$'s under the convexity of paths of order three. \emph{Discrete Appl.\ Math.} \textbf{192} (2015), 28--39.
\href{https://doi.org/10.1016/j.dam.2014.05.005}{doi:10.1016/j.dam.2014.05.005}

\bibitem{cap-2019} M.R.~Cappelle, E.M.M.~Coelho, H.~Coelho, B.R.~Silva, F.~Protti, U.E.~Souza. $P_3$-hull number of graphs with diameter two. \emph{Electron.\ Notes Theor.\ Comput.\ Sci.} \textbf{346} (2019), 309--320.
\href{https://doi.org/10.1016/j.entcs.2019.08.028}{doi:10.1016/j.entcs.2019.08.028}

\bibitem{cen-2010} C.C.~Centeno, S.~Dantas, M.C.~Dourado, D.~Rautenbach, J.L.~Szwarcfiter. Convex partitions of graphs induced by paths of order three. \emph{Discrete Math.\ Theor.\ Comput.\ Sci.} \textbf{12} (2010), 175--184.
\href{https://doi.org/10.46298/dmtcs.502}{doi:10.46298/dmtcs.502}

\bibitem{cha-1979} J.~Chalupa, P.L.~Leath, G.R.~Reich. Bootstrap percolation on a Bethe latice. \emph{J. Phys. C} \textbf{12} (1979), 31--35.
\href{https://doi.org/10.1088/0022-3719/12/1/008}{doi:10.1088/0022-3719/12/1/008}

\bibitem{coe-2019} E.M.M.~Coelho, H.~Coelho, J.R.~Nascimento, J.L.~Szwarcfiter. On the $P_3$-hull number of some products of graphs. \emph{Discrete Appl.\ Math.} \textbf{253} (2019), 2--13.
\href{https://doi.org/10.1016/j.dam.2018.04.024}{doi:10.1016/j.dam.2018.04.024}

\bibitem{coe-2015} E.M.M.~Coelho, M.C.~Dourado, R.M.~Sampaio. Inapproximability results for graph convexity parameters. \emph{Theoret.\ Comput.\ Sci.} \textbf{600} (2015), 49--58.
\href{https://doi.org/10.1016/j.tcs.2015.06.059}{doi:10.1016/j.tcs.2015.06.059}

\bibitem{gon-2021} L.N.~Gonz\'alez, L.N.~Grippo, M.D.~Safe. Formulas in connection with parameters related to convexity of paths on three vertices: caterpillars and unit interval graphs. \emph{Australas.\ J.\ Combin.} \textbf{79} (2021), 401--423.
\href{https://doi.org/10.48550/arXiv.2002.10598}{doi:10.48550/arXiv.2002.10598}

\bibitem{gri-2021} L.N.~Grippo, A.~Pastine, P.~Torres, M.~Valencia-Pabon, J.C.~Vera. On the $P_3$-hull number of Kneser graphs. \emph{Electron.\ J.\ Combin.} \textbf{28} (2021), \#P3.32.
\href{https://doi.org/10.37236/9903}{doi:10.37236/9903} 

\bibitem{ham-book} R.~Hammack, W.~Imrich, S.~Klav\v{z}ar. \emph{Handbook of Product Graphs}. CRC Press, Boca Raton, FL, 2011. 

\bibitem{morrison} M.~Morrison, J.~Noel. Extremal bounds for bootstrap percolation in the hypercube. \emph{J.\ Comb.\ Theory Series A} \textbf{156} (2018), 61--84.
\href{https://doi.org/10.1016/j.jcta.2017.11.018}{doi:10.1016/j.jcta.2017.11.018}

\bibitem{prz-2012} M.~Przykucki. Maximal percolation time in hypercubes under $2$-bootstrap percolation. \emph{Electron.\ J.\ Combin.} \textbf{19} (2012), \#P41.
\href{https://doi.org/10.37236/2412}{doi:10.37236/2412}

\bibitem{prz-2020} M.~Przykucki, T.~Shelton. Smallest percolating sets in bootstrap percolation on grids. \emph{Electron.\ J.\ Combin.} \textbf{27} (2020), \#P4.34.
\href{https://doi.org/10.37236/9582}{doi:10.37236/9582}

\bibitem{shitov} Y.~Shitov. Counterexamples to Hedetniemi's conjecture. \emph{Ann.\ Math.} \textbf{190} (2019), 663--667.
\href{https://doi.org/10.4007/annals.2019.190.2.6}{doi:10.4007/annals.2019.190.2.6}

\bibitem{tar} C.~Tardif. Hedetniemi's conjecture, 40 years later. \emph{Graph Theory Notes N. Y.} \textbf{54} (2008), 46--57. 

\bibitem{west} D.B.~West. \emph{Introduction to Graph Theory (Second Edition)}. Prentice Hall, USA, 2001.

\bibitem{zhu} X.~Zhu. A survey of Hedetniemi's conjecture. \emph{Taiwanese J. Math.} \textbf{2} (1998), 1--24.
\href{https://doi.org/10.1016/j.cosrev.2025.100794}{doi:10.1016/j.cosrev.2025.100794}


\end{thebibliography}
\end{document}